\documentclass[12]{article}

\usepackage{amsthm, amssymb, amsmath}
\usepackage{amscd}
\usepackage{epsf}
\usepackage{latexsym}
\usepackage{verbatim}
\usepackage[dvips]{graphicx}

\newtheorem{thm}{Theorem}
\newtheorem{prop}[thm]{Proposition}
\newtheorem{cor}[thm]{Corollary}
\newtheorem{lem}[thm]{Lemma}
\newtheorem{remark}[thm]{Remark}

\theoremstyle{plain}

\swapnumbers
\numberwithin{thm}{section}

\newcommand{\Z}{\mathbb{Z}}
\newcommand{\N}{\mathbb{N}}

\newcommand{\If}{\textrm{if}}

\newcommand{\No}{\mathbb{N}_0}
\newcommand{\half}{\frac{1}{2}}

\newcommand{\Rm}{\mathbb{R}}

\title{Ergodic Properties of a Class of Discrete Abelian Group Extensions of Rank-One Transformations}

\author{Chris Dodd
 \thanks {Department of Mathematics, Massachusetts Institute of Technology, MA, 02139, USA.
    cdodd@math.mit.edu}
 \and Phakawa Jeasakul
 \thanks{Economics Department, University of California, Berkeley, CA 94720, USA.
 phakawa@econ.Berkeley.edu}
  \and Anne Jirapattanakul
 \thanks{1022 International Affairs Building, 
  Columbia University,
  420 West 118th Street, New York, NY 10027, USA.
pj2133@columbia.edu}
\and
{Daniel M. Kane} \thanks {Department of Mathematics, Harvard University, Cambridge, MA 02139, USA.  dankane@math.harvard.edu}
  \and Becky  Robinson    \thanks{Williams College, Williamstown, MA 01267, USA.
    05err@williams.edu}
  \and Noah Stein
 \thanks {Laboratory for Information and Decision Systems, Massachusetts Institute of Technology, Cambridge, MA 02139, USA.
    nstein@mit.edu}
    \and Cesar E.  Silva
     \thanks{Department of Mathematics,
    Williams College, Williamstown, MA 01267, USA.
    csilva@williams.edu}
  }

\begin{document}

\maketitle

\begin{abstract}
We define a class of discrete abelian group extensions of
rank-one transformations and establish necessary and
sufficient conditions for these extensions to be power weakly
mixing.  We  show that all members of this class
are multiply recurrent.   We then  study conditions sufficient for showing that
cartesian products of transformations are conservative for a class
of invertible infinite measure-preserving transformations and
provide examples of these transformations.
\end{abstract}

\section{Introduction}

Group extensions of measure-preserving dynamical systems have received much attention in the literature.  In most of the works the group has been assumed to be compact, and if the base transformation
is finite measure-preserving then the extension is finite measure-preserving.  A
question  that has been studied in this context  is conditions under which dynamical properties of the
base transformation (such as weak mixing or mixing) lift to the group extension; the reader may refer to e.g.  \cite{Ru}, \cite{Ro} and the references in these works.   In this article we consider  extensions of a class
of rank-one transformations by countable discrete abelian groups.  While the base transformation is restricted to be a rank-one transformation we allow the group to possibly be infinite.   We establish a simple condition   that is equivalent to the ergodicity of the  extensions, and another condition that is equivalent to power weak mixing of the extensions.  Power weak mixing is equivalent to weak mixing for finite measure-preserving transformations,  but it is a stronger property in the case  of infinite measure-preserving transformations.
We show that the extension is power weakly mixing if it is totally
ergodic.  We also show that our group extensions are multiply recurrent, and give several applications showing ergodicity or (power) weak mixing for certain extensions in both the finite and infinite  measure-preserving cases.  In the later sections we consider the question of the conservativity of products of powers of infinite measure-preserving transformations,  and apply our results to staircase transformations.

Let $(X,\mathcal B,\mu)$ be a measure space isomorphic to a finite or infinite interval in
$\mathbb R$ with Lebesgue measure $\mu$ (when the interval is finite we  assume $\mu$ has been normalized to be a probability measure).  Let $T: X\to X$ be an invertible measure-preserving transformation.
The transformation $T$ is {\bf conservative} if for any set $A$ of
positive measure, there exists an integer $i>0$ such that
$\mu(T^{-i} A \bigcap A) > 0$.  $T$ is {\bf ergodic } if for any pair of sets set $A$ and $B$ of
positive measure, there exists an  integer $i\geq 0$ such that
$\mu(T^{-i} A \bigcap B) > 0$.  (As our transformations are invertible and defined on nonatomic spaces, ergodicity
implies conservativity.)  Let
$T^{\otimes d}$ denote  the cartesian product of $d>0$ copies of $T$.
We say that   $T$ has {\bf infinite conservative
index} if $T^{\otimes d}$ is conservative with respect to $\mu^d$ for all $d>0$ (where $\mu^d$ denotes $d$-dimensional product of $\mu$);
$T$ has {\bf infinite ergodic index} if for all $d>0$, $T^{\otimes d}$ is ergodic with respect to $\mu^d$.
A transformation $T$ has {\bf power conservative
index} if for all sequences of positive integers $k_1,k_2,\ldots k_d, T^{k_1} \times T^{k_2} \times\ldots\times T^{k_d} : X^{\otimes d}
\rightarrow X^{\otimes d}$ is conservative; $T$ is said to be {\bf power weakly mixing}
if for all   nonzero $k_1,\ldots, k_d$, $ T^{k_1} \times T^{k_2} \times\ldots\times T^{k_d} $ is ergodic.

Power weak mixing is clearly equivalent to weak mixing for finite measure-preserving transformations, but it is a stronger property in the case  of infinite measure-preserving transformations \cite{AFS01}.   In fact, there exists a transformation  $T_1$ such that $T_1$ has infinite ergodic index but  $T_1 \times T_1^2$ is not conservative, hence not ergodic (\cite{AFS01}).

In Section~\ref{S:prelim} we define, for each countable discrete abelian group $G$, a class of
measure-preserving transformations.   When $G$ is an infinite group the transformation
is infinite measure-preserving.  As the last example in  Section~\ref{S:examples}
shows, these contain group extensions of rank-one transformations.  In Theorem~\ref{pwmthm} we give
necessary and sufficient conditions for our construction to be power weakly mixing.
When $G$ is a finite group, the transformation is finite measure-preserving and our
theorem gives equivalent conditions for weak mixing.

We also show that our group extensions are multiply recurrent.  A transformation $T$ is said to be
{\bf $d$-recurrent}  if for all sets of positive measure $A$ there exists an integer $n>0$ such that
$\mu(A\cap T^n(A)\cap\cdots\cap T^{nd}(A))>0$. $T$ is said to be {\bf multiply recurrent} if it is $d$-recurrent for all integers $d>0$.
  As is well-known, Furstenberg showed that every finite measure-preserving transformation is multiply recurrent \cite{hF81}, but it is now known that  infinite measure-preserving transformations need not be multiply recurrent \cite{EHH98}, \cite{AN00}, even when they are power weakly mixing \cite{G03}.  However, it was shown recently that compact group \cite{inoue} and $\sigma$-finite \cite{Mey}  extensions of
multiply recurrent infinite measure-preserving transformation are multiply recurrent, .  This need not be the case for extensions by non-compact groups as already observed in \cite{inoue},  but we obtain multiple recurrence for  our class of (non $\sigma$-finite) extensions.  In particular, it follows  that for each countable discrete abelian group there is
a multiply recurrent extension.

In Section~\ref{S:pci} we introduce a condition for rank-one transformations that implies power conservative index,  and use it show show that some infinite measure-preserving staircases have power conservative index.

\bigskip
\noindent {\bf Acknowledgments.}  This paper is based on research in
the Ergodic Theory group of the 2004 SMALL Undergraduate Summer
Research Project at Williams College, with Silva as faculty advisor.
Support for the project was provided by a National Science Foundation
REU Grant and the Bronfman Science Center of Williams College.


\section{Construction of the Transformations}\label{S:prelim}

  Fix a countable discrete abelian group $G$.  We will
construct transformations that are $G$ extensions of rank-one
transformations produced by a standard cutting and stacking
procedure.  Let $\Gamma$ be the set of  all
elements that are of the form
\[
(\gamma_e,s_{e,0},...,s_{e,\gamma_e-1},g_{e,0},...,g_{e,\gamma_e-1})\]
where  $\gamma_e>1$ is  natural number and the remaining entries are an 
element of $\N_0^{\gamma_e}\times G^{\gamma_e}$.  For clarity, we sometimes
write the subscript $s_{e,0}$ as $s(e,0)$, etc.  We think of $\Gamma$ as the set of possible operations to go
from one generation to the next.  $\gamma$ is the number of pieces
that we cut each level into, $s_{e,i}$ describes the numbers of
spaces added, and $g_{e_i}$ describe how the $G$-component of the
column changes.  Let
\[
F:\No \rightarrow \Gamma
\] be
a function. We think of $F$ as the map from generation
numbers to what operation is performed in that generation.  We
require that $F$ have the property that for any natural numbers $n$
and $d$, there are infinitely many natural numbers $m$ so that
$F(n+i)=F(m+i)$ for all $0\leq i < d$.  In other words any
sequence that appears in $F$ does so infinitely often. Let \[F(n) =
(\gamma_n,s(n,0),...,s(n,\gamma_n -1),g(n,0),...,g(n,\gamma_n-1)).\]
Given $F$, we define a (at most rank-$|G|$) transformation $T$ as follows:

A {\bf{ column}} consists of a finite (ordered) sequence of intervals of the same length, called the  {\bf {levels}}
of the column; the number of levels is the \emph{height} of the column.  We begin with generation-$0$ columns $C_{0,g}$ for $g\in G$, each
consisting of an interval  of mass $1$.  To obtain the generation-$(N+1)$ columns from the
generation-$N$ columns,  first write each generation-$N$ column $C_{N,g}$ ($g\in G$) as
\[C_{N,g}=(I_{N,g}^{(0)},I_{N,g}^{(1)},...,I_{N,g}^{(h_{N}-1)}),\]
where we think of $h_N$ as the height of the column.
 Next,
 cut each level or interval  $I_{N,g}^{(i)}$ into $\gamma_N$ equal mass subintervals
\[I_{N,g,0}^{(i)}, \ldots,I_{N,g,\gamma_N-1}^{(i)},\]   and set
\begin{align*}
C_{N+1,g} = &
(I_{N,g+g(N,0),0}^{(0)},\ldots,I_{N,g+g(N,0),0}^{(h_{N}-1)},S_{N,g,0}^{(0)},\ldots,S_{N,g,0}^{(s(N,0)-1)},\\
& I_{N,g+g(N,1),1}^{(0)},\ldots,I_{N,g+g(N,1),1}^{(h_{N}-1)},
S_{N,g,1}^{(0)},\ldots,S_{N,g,1}^{(s(N,1)-1)},\ldots,\\
&
I_{N,g+g(N,\gamma_N-1),\gamma_N-1}^{(0)},\ldots,I_{N,g+g(N,\gamma_N-1),\gamma_N-1}^{(h_{N}-1)}
,\\
&\quad\quad\quad \quad\quad\quad\quad\quad\quad \quad\quad\quad S_{N,g,\gamma_N-1}^{(0)},\ldots,S_{N,g,\gamma_N-1}^{(s(N,\gamma_N-1)-1)})
\end{align*}
where each $S_{N,g,i}^{(j)}$ is a spacer level, i.e., a new subinterval of the same
length as any of the subintervals in its column.   The resulting  transformation is defined
on the intervals of each column by sending that interval by translation to the
interval above it if there is one.  In the limit, the lengths of the intervals in each column
converges to zero, so  the transformation is defined in the union of all the levels.
We thus obtain a transformation  $T$ that is measure preserving.  Furthermore, one can
arrange the subintervals in each column so that $T$ is defined on
a finite or infinite subinterval of $\mathbb R$.



We prove the following theorems:

\begin{thm} \label{multrecthm}
For all such $F$, $T$ is multiply recurrent.
\end{thm}

\begin{thm} \label{pwmthm}
$T$ is power weakly mixing if and only if the following conditions
are both satisfied,
\begin{enumerate}
\item $\{ g(N,i)-g(N,0): N\in \No, 0\leq i \leq \gamma_N-1 \}$ generate $G$
\item For all $N$, $(1,0)$ is in the integer span of
\begin{align*}
&\{ (s(N,i)+h_N,g(N,i+1)-g(N,i)) : 0\leq i \leq \gamma_N -2\}
\cup \\
&\{((s(M+1,i)+s(M,\gamma_M-1)-s(M,0)),\\
&(g(M+1,i+1)-g(M+1,i)+2g(M,0)-g(M,\gamma_M-1)-g(M,1)))\\
&: M\in \No, 0\leq i \leq \gamma_{M+1}-2\}
\end{align*}
in $\Z \times G$.
\end{enumerate}
\end{thm}

The first condition essentially states that it is possible to get from any column to any other column.  The $\Z \times G$ that appears in the second condition should be thought of as a group acting on our space with $G$ acting by changing column index, and $1\in\Z$ acting as $T$.  Let us call the terms in the second condition
$$
t_{N,i} = (s(N,i)+h_N,g(N,i+1)-g(N,i))
$$
and
\begin{align*}
c_{M,i} = (&(s(M+1,i)+s(M,\gamma_M-1)-s(M,0)),\\
&(g(M+1,i+1)-g(M+1,i)+2g(M,0)-g(M,\gamma_M-1)-g(M,1))).
\end{align*}
They each represent distances in this action between copies of columns as will be discussed later.  The condition that $(1,0)$ be in their span essentially says that we have the control to shift things by $T$.

\section{Some Machinery Involving copies of Columns}

If $I$ is a level of a generation-$n$ column,
 $n > 1$, we say that a level  $K$ in a generation-$(n+m)$ column is 
a {\bf{copy}} of $I$ if $K$ corresponds to a subset of level $I$.
We define a copy of a column $C$, in some column of later
generation, to be a union of consecutive levels that are, in order,
copies of the levels of $C$.
We would like to be able to index the copies of generation-$N$
columns in a particular generation-$(N+M)$ column.  If $C$ is a copy
of $C_{N+1,g}$, then we let $P_i (C)$ be the copy of
$C_{N,g+g(N,i)}$ contained in $C$. In particular,  for $0\leq i \leq \gamma_N-1$ $P_i(C)$ is
the $i^{th}$ copy of a generation-$N$ column contained  in $C_{N,g}$.  Let
\[
P_{N,g}[a_0,a_1,...,a_n] = P_{a_0}(P_{a_1}(\ldots P_{a_n}(
C_{N,g})\ldots )),
\]
where $C_{N,g}$ is thought of as a copy of itself.

Notice that the $P_{N,g}[a_0,\ldots,a_{n-1}]$ index all of the copies of generation-$(N+n)$ columns in $C_{N,g}$.  Their relative positions are given by the radix ordering on the $a_i$ with $a_0$ being the most significant.

\begin{lem} \label{copycolor}
$P_{N+n,g}[a_0,a_1,...,a_{n-1}]$ is a copy of
$C_{N,g+\sum_{i=0}^{n-1} g(N+i,a_i)}$.
\end{lem}
\begin{proof}
We proceed by induction on $n$.  The $n=0$ case is trivial. Assuming
that our statement holds for $n-1$, we have, letting $C'_{M,h}$
denote a copy of $C_{M,h}$ for any $M\in\No$ and $h\in G$, that
\begin{align*}
P_{N+n,g}[a_0,a_1,...,a_{n-1}] = & P_{a_0}
\left(C'_{N+1,g+\sum_{i=1}^{n-1} g(N+i,a_i)}\right)\\
= & C'_{N,g+\sum_{i=0}^{n-1} g(N+i,a_i)}.
\end{align*}
This completes our inductive step and proves our Lemma.
\end{proof}

\begin{lem} \label{totaldistancelem}
$T^k(P_{N+n,g}[a_0,\ldots,a_{n-1}])=P_{N+n,g}[b_0,\ldots,b_{n-1}]$
where
$$
k = \sum_{i=0}^{n-1}\left(h_{N+i}(b_i-a_i) + \sum_{j=0}^{b_i-1}
s(N+i,j) - \sum_{j=0}^{a_i-1} s(N+i,j) \right).
$$
\end{lem}
\begin{proof}
We proceed by induction on $\sum_{i=0}^{n-1} |a_i-b_i|$.  The
statement is clearly true when this is 0.  Otherwise, assuming our
hypothesis for smaller values of $\sum_{i=0}^{n-1} |a_i-b_i|$.
Without loss of generality we may assume that $b_i>a_i$.  Then we
have that
\begin{align*}
P_{N+n,g}[b_0,\ldots,b_{n-1}] &
=T^{h_{N+i}+s_{b_i-1}}(P_{N+n,g}[b_0,\ldots b_{i-1},b_i -1,
b_{i+1},\ldots,b_{n-1}])\\
& = T^k(P_{N+n,g}[a_0,\ldots,a_{n-1}]).
\end{align*}
This completes our inductive step and proves the Lemma.
\end{proof}

\begin{lem} \label{consecdistancelem}
$T^k
(P_{N+n,g}[\gamma_N-1,\gamma_{N+1}-1,\ldots,\gamma_{N+m}-1,a_{m+1},\ldots,a_{n-1}])
=P_{N+n,g}[0,\ldots,0,a_{m+1}+1,a_{m+2},\ldots,a_{n-1}]$ where
$a_{m+1}\leq \gamma_{N+m+1}-2$ and
$$
k=h_N+\sum_{i=0}^{m} s(N+i,\gamma_{N+i}-1) + s(N+m+1,a_{m+1}).
$$
\end{lem}
\begin{proof}
By Lemma \ref{totaldistancelem} we have that
$$
k = \sum_{i=0}^{m} \left( -(\gamma_{N+1}-1)h_{N+i} -
\sum_{j=0}^{\gamma_{N+i}-2} s(N+i,j)\right) +
h_{N+m+1}+s(N+m+1,a_{m+1}).
$$
Using the fact that $h_{n+1} = \gamma_n h_n +
\sum_{j=0}^{\gamma_n-1}s(n,j)$, we have that
\begin{align*}
k &=h_{N+m+1}+s(N+m+1,a_{m+1}+\sum_{i=0}^m h_{N+1} +
s(N+i,\gamma_{N+i}-1) - h_{N+i+1}\\
& = h_N+\sum_{i=0}^{m}
s(N+i,\gamma_{N+i}-1) + s(N+m+1,a_{m+1}).
\end{align*}
Thus proving our Lemma.
\end{proof}

$t_{N,i}$ represents the change in location of the copy when the index corresponding to generation-$N$ is changed from $i$ to $i+1$.  $c_{M,i}$ represents the change in location of the copy corresponding to the pair of indices corresponding to generations $M$ and $M+1$ changing from $(\gamma_M-1,i)$ to $(0,i+1)$.

\section{Necessity of the Conditions}

\begin{lem} \label{erodiclem}
$T$ is ergodic only if condition 1 is satisfied.
\end{lem}
\begin{proof}
Suppose that $\{ g(N,i) : N\in\No, 0\leq i \leq \gamma_N -1\}$
generate $H \subsetneq G$.  Let $g_1,g_2\in G$ be in different
cosets of $H$.  Consider $A=C_{0,g_1},B=C_{0,g_2}$.  Assume for sake 
of contradiction, that for some $n$, $\mu(T^n(A)\cap B) > 0$.  This
would imply that there is some column that contains both a copy of
$C_{0,g_1}$ and a copy of $C_{0,g_2}$.  Suppose these copies are
$P_{N,g_3}[a_0,\ldots,a_{n-1}],P_{N,g_3}[b_0,\ldots,b_{n-1}]$.  Then
by Lemma \ref{copycolor} we have that $g_1 = g_3 + \sum_{i=0}^{n-1}
g(i,a_i)$ and $g_2 = g_3 + \sum_{i=0}^{n-1} g(i,b_i)$.  Hence we
have that $g_2-g_1 = \sum_{i=0}^{n-1} g(i,b_i)-g(i,a_i) \in H$, but
this is not the case.  Hence $T$ is not ergodic.
\end{proof}

\begin{remark}
This condition is actually sufficient for $T$ being ergodic.  This
fact will follow from Lemma \ref{samecollem}.
\end{remark}

\begin{lem} \label{totalergodicitylem}
$T$ is totally ergodic only if condition 2 is satisfied.
\end{lem}
\begin{proof}
Pick an $N$ for which the condition does not hold.
Let $t_{N,i}= t_i$.
Let
\begin{align*}
H  & = \textrm{span} (\{t_i : 0\leq i \leq \gamma_N -2\} \cup
\{c_{M,i} : M\in \No, 0\leq i \leq \gamma_{M+1}-2\})
\end{align*}
Let $H\cap \left(\Z\times \{ 0 \}\right) \subset \Z (D , 0)$ for $D>1$.  We will
prove that there is no integer, $n$, so that
$\mu(T^{nD}(I_{N,0}^{(0)})\cap I_{N,0}^{(1)})>0.$  Suppose for sake
of contradiction that this is not the case.  Then there must exist
copies $C_1,C_2$ of $C_{N,0}$ that are in the same column, and with
$T^{nD-1}(C_1) = C_2$.  Suppose that these copies are
$P_{N+l,g}[a_0,\ldots,a_{l-1}],P_{N+l,g}[b_0,\ldots,b_{l-1}]$.  For
two copies of generation-$N$ columns in $C_{N+l,g}$, $\alpha,\beta$
define $\Delta(\beta,\alpha)=(k,h)$ where $T^k(\alpha)=\beta$, and
$\alpha$ and $\beta$ are copies of $C_{N,g'}$ and $C_{N,g'+h}$
respectively.  Notice that
$\Delta(\gamma,\beta)+\Delta(\beta,\alpha) = \Delta(\gamma,\alpha)$.
Notice also that $\Delta(C_2,C_1)=(nD-1,0)$.  Lastly, notice that
$C_1$ and $C_2$ are connected by some chain of copies where each
pair of consecutive copies are of the form given in Lemma
\ref{consecdistancelem}.  We have by Lemmas \ref{copycolor} and
\ref{consecdistancelem} that
\begin{align*}
\Delta( & P_{N+l,g}[0,\ldots,0,a_{m+1}+1,a_{m+2},\ldots,a_{l-1}],\\
&
P_{N+l,g}[\gamma_N-1,\gamma_{N+1}-1,\ldots,\gamma_{N+m}-1,a_{m+1},\ldots,a_{l-1}])
= \\
( h_N & + s(N+m+1,a_{m+1}) +\sum_{i=0}^{m} s(N+i,\gamma_{N+i}-1), \\
& g(N+m+1,a_{m+1}+1)-g(N+m+1,a_{m+1})-\sum_{i=0}^m
g(N+i,\gamma_{N+i}-1)) = \\
& t_0 + \sum_{i=0}^{m-1} c_{N+i,0} + c_{N+m,a_{m+1}} \in H.
\end{align*}
Combining these facts with the fact that $H$ is additively closed,
we have that $\Delta(C_2,C_1)\in H$.  But by assumption,
$(nD-1,0)\notin H$. This is a contradiction.  Hence $T^D$ is not
ergodic, proving our Lemma.

\end{proof}

\section{Sufficiency of Conditions}

\begin{lem} \label{technicallem}
If condition 2 is satisfied then for any $N\in \No$, there 
exists some
 $D\neq 0$  so that  $(D,0,0)$ is in
the integer span of
\begin{align*}
&\{ (s(N,i),g(N,i+1)-g(N,i),1) : 0\leq i \leq \gamma_N -2\}
\cup \\
&\{(c_{M,i},0): M\in \No, 0\leq i \leq \gamma_{M+1}-2\}
\end{align*}
in $\Z\times G\times\Z$ for any $N\in \No$.
\end{lem}
\begin{proof}
Let the intersection of the integer span of this with $\Z\times \{ 0
\} \times \Z$ be $H$.  Consider the homomorphism, $\phi_M:\Z\times
G\times\Z \rightarrow \Z \times G$ defined by $\phi_M (a,b,c)=(a+h_M
c,b)$.  Notice that $\phi$ sends $\Z\times\{ 0\}\times\Z$ to
$\Z\times\{ 0\}$.  Hence for all $M$ with $F(M)=F(N)$, we have that
$(1,0)\in \phi_M(H)$.  Suppose for sake of contradiction that $H\cap
(\Z \times \{0\} \times \{ 0\}) = \{(0,0,0)\}$.  Then $H$ must have
infinite index in $\Z \times \{ 0 \} \times \Z$.  So $H=\Z(a,0,b)$
for some $a$ and $b$. But then, $\phi_M(H) = \Z(a+b h_M,0)$.  But
this implies that for infinitely many values of $h_M$ that $a+b h_M
= \pm 1$. This implies that $b=0$ and $a = \pm 1$.  But then
$(1,0,0) \in H \cap (\Z \times \{0\} \times \{ 0\}).$  This is a
contradiction.
\end{proof}

The remainder of this section is devoted to proving that
$T^{k_1}\times\cdots\times T^{k_m}$ is ergodic.  We will let $A$ and
$B$ be arbitrary sets of positive measure in the $m$-fold product of
the space on which $T$ is defined.   Given an interval $I$,
a measurable set $A$ and $\delta>0$, we say that $I$ is
more than $(1-\delta)$-full of $A$ if $\mu(A\cap I)>(1-\delta)\mu(I)$;
a similar notion is defined for product sets.

We use a standard technique from measure theory, sometimes
called double approximation, which states that if we have sets
$I$ and $J$ (generally products of levels) more than half-full of
$A$ and $B$, respectively, then if we have several generations of
partitions of $I$ and $J$ into equal numbers of subsets of equal
measure, each generation a refinement of the last (generally, the
partitions consist of all products copies of intervals form $I$ and
$J$ in some later generation), with bijections between the $N^{th}$
generation subsets of $I$ and the $N^{th}$ generation subsets of
$J$, and if furthermore these subsets form a basis of the topologies
on $I$ and $J$, then for any $\epsilon>0$, there exist corresponding
subsets of $I$ and $J$ of some generation that are
$(1-\epsilon)$-full of $A$ and $B$ respectively (see e.g. \cite[6.5.4]{Si08}).

\begin{lem} \label{samecollem}
If condition 1 holds, then there exist an integer $N$ and levels
$I_1,\ldots,I_m$, $J_1,\ldots,J_m$ from generation-$N$ columns so
that $I_i$ and $J_i$ are in the same column and so that
$I_1\times\cdots\times I_m$ and $J_1\times\cdots\times J_m$ are more
than $\left( \half \right)$-full of $A$ and $B$ respectively.
\end{lem}
\begin{proof}
We may find an integer $N_1$ and generation-$N_1$ columns
$I'_1,\ldots,I'_m, $ $J'_1,\ldots,J'_m$ so that $I'_1\times\cdots\times
I'_m$ and $J'_1\times\cdots\times J'_m$ are more than $\left( \half
\right)$-full of $A$ and $B$ respectively.  Let $I'_i$ be in
$C_{N_1,h_{1,i}}$ and let $J'_i$ be in $C_{N_1,h_{2,i}}$.  By
condition 1, we may write
$$
h_{1,i}+\sum_{j=1}^{r_{1,i}}g(e_{1,i,j},l_{1,i,j})=h_{2,i}+\sum_{j=1}^{r_{2,i}}g(e_{2,i,j},l_{2,i,j})
$$
for some values of $r\in\No$, $e\in S$ and $l\in\No$.  Since $F$
attains all values in $S$ infinitely often, we may find some
sequence of consecutive integers, $a,a+1,\ldots,a+b$, so that for
each $t \in \{1,2\}$, $1\leq i \leq m$, and $1\leq j \leq r_{t,i}$
there is a distinct $0\leq \alpha_{t,i,j} \leq b$ so that
$F(a+\alpha_{t,i,j})=e_{t,i,j}$.  Using double approximation, we may
find an integer $N_2$ so that $F(N_2+i)=F(a+i)$ for all $0\leq i
\leq b$ and generation $N_2$ levels
$I''_1,\ldots,I''_m,J''_1,\ldots,J''_m$ so that
$I''_1\times\cdots\times I''_m$ and $J''_1\times\cdots\times J''_m$
are more than $\left( \frac{1}{2\prod_{j=a}^{a+b} \gamma_j^m}
\right)$
-full of $A$ and $B$ respectively.  Furthermore, we can
ensure that if $I''_i$ and $J''_i$ are in columns $C_{N_2,h'_{1,i}}$
and $C_{N_2,h'_{2,i}}$ respectively, that
$h'_{1,i}-h'_{2,i}=h_{1,i}-h_{2,i}$.  Then if we let $N=N_2+b+1$ and
let $I_i$ be the copy of $I''_i$ in
$$
P_{N,h'_{1,i}-\sum_{j=1}^{r_{2,i}}g(e_{2,i,j},l_{2,i,j})}[d_0\ldots
d_b]
$$
where
$$
d_p = \begin{cases} l_{2,i,j} \ \ & \If \ p=\alpha_{2,i,j}\\
0 \ \ & \If \ p\neq \alpha_{2,i,j} \forall j
\end{cases}
$$
and $J_i$ be the copy of $J''_i$ in
$$
P_{N,h'_{2,i}-\sum_{j=1}^{r_{1,i}}g(e_{1,i,j},l_{1,i,j})}[d'_0\ldots
d'_b]
$$
where
$$
d'_p = \begin{cases} l_{1,i,j} \ \ & \If \ p=\alpha_{1,i,j}\\
0 \ \ & \If \ p\neq \alpha_{1,i,j} \forall j
\end{cases}.
$$
These are copies of the correct columns by Lemma \ref{copycolor}.
They are clearly in the same column.  Furthermore we have that
$I_1\times\cdots\times I_m$ and $J_1\times\cdots\times J_m$ are more
than $\left( \half \right)$-full of $A$ and $B$ respectively,
proving our Lemma.
\end{proof}

\begin{lem}\label{multdislem}
If conditions 1 and 2 hold, there exist levels
$I_1,\ldots,I_m,J_1,\ldots,J_m$ that satisfy the conditions from
Lemma \ref{samecollem} with the additional property that some power
of $T^D$ sends $I_i$ to $J_i$, where $D$ satisfies the
statement of Lemma \ref{technicallem} for some $N_0$
\end{lem}
\begin{proof}
The proof follows the same lines as that of Lemma \ref{samecollem}.
We start with the levels given to us by Lemma \ref{samecollem} and
then use double approximation to get the levels that we need.  
Suppose that $I_i$ and $J_i$ are separated by $T^{r_i}$.  Using double approximation we know that for and $\epsilon>0$ and all sufficiently large generation numbers $N$, we can find generation-$N$ copies $I_i'$ and $J_i'$ of $I_i$ and $J_i$ respectively, so that $T^{r_i}(I_i') = J_i'$ and so that $I_1'\times\ldots\times I_m'$ and $J_1'\times\ldots\times J_m'$ are at least $\left(1-\epsilon^m/2\right)$-full of $A$ and $B$ respectively.  What we wish to show is that for some $\epsilon>0$, and for arbitrarily large generation numbers $N$, given any such intervals $I_i'$ and $J_i'$, that we can find copies $I_i''$ and $J_i''$ of these in generation-$(N+n)$, that are of size at least $\epsilon$ that of the original, and so that $I_i''$ and $J_i''$ are in the same column, separated by a power of $T^D$.  The result would then follow since $I_1''\times\ldots\times I_m''$ and $J_1''\times\ldots\times J_m''$ would be at least $\left( \frac{1}{2}\right)$-full of $A$ and $B$ respectively.

For the above to work, we need only show that for any separation $r=r_i$, that for some sufficiently small $\epsilon>0$ and sufficiently large generation $M$, that we can find two copies of $C_{M,g}$ that are of size at least $\epsilon$ that of the original, are in the same column, and are separated by a power of $T$ congruent to $r$ modulo $D$.  This allows us to produce the necessary copies of $I'_i$ and $J_i'$ for each $i$.


For each congruence class, $c$ modulo $D$ such that for infinitely many $N$, $h_N$ is in $c$ and $F(N)=F(N_0)$,  we can, by
condition 2, find some integer combination of 
$t_{N,i}$ and $c_{M,i}$ That add up to $(r,0)$ in $(\Z/D\Z )\times G$.  Suppose that the sum
of the absolute values of the multiples of terms of the form
$t_{N,i}$ needed is at most $X$.  For $a,b\in S$
suppose that the sum of the absolute values of multiplies terms of the form
$c_{M,i}$ where $F(M)=a$ and $F(M+1)=b$ needed is at most $Y_{a,b}$.  Find a
string of consecutive integers, $I$, so that on this string the
following hold:

There are at least $XD$ values $n\in I$ so that $F(n)=e$.

There are a number of non-overlapping pairs of consecutive integers
in $I$ which do not intersect any of the $n$ used in the previous
condition, so that for at least $2DY_{a,b}$ of these pairs, $F$
evaluated at these values yields $a$ and $b$ in that order.

Then for any interval of sufficiently large numbers, $I'$ on which
$F$ agrees with the values it takes on $I$, we can find one of these
congruence classes, $c$ for which there are at least $X$ values
$n\in I'$ for which $F(n)=e$ and $h_n \equiv c \pmod{D}$.

For each $a$ and $b$ we can find at least $2Y_{a,b}$ pairs of
consecutive integers $n,n+1\in I'$ so that $F(n)=a$, $F(n+1)=b$ and
$h_n$ has the same value modulo $D$ for all of these pairs.

Now if $M'$ is the smallest value in $I'$, we can construct two copies
of $C_{M',g}$  whose size is at least
$\prod_{j\in I} \frac{1}{\gamma_j}$ of the original.  Suppose that
$$
 \sum_{i=0}^{n-2} \alpha_i (s(N,i)+c,g(N,i+1)-g(N,i)) +
 \sum_{M,i} \beta_{M,i} c_{M,i} = (r,0)
$$
In $(\Z/D\Z )\times G$.  Then we consider copies of the form
$$
P_{M'+k+1,g}[d_0\ldots d_k],P_{M'+k+1,g}[d'_0\ldots d'_k]
$$
where $M'+k$ is the largest value in $I'$.  We define the
$d_i$ and $d'_i$ as follows:

There are $\alpha_i$ values $n\in I'$ for which $F(n)=e$ and $h_n
\equiv c \pmod{D}$ where $d'_{n-M'} = i+1$ and $d_{n-M'} = i$ (if
$\alpha_i$ is negative, we reverse the values and do it $|\alpha_i|$
times).

There are $\beta_{M,i}$ values $n\in I'$ where $F(n)=F(M)$,
$F(n+1)=F(M+1)$, $d'_{n-M'}=0$, $d_{n-M'}=\gamma_M-1$,
$d_{n-M'+1}=i$, $d'_{n-M'+1}=i+1$, and the same number of such
values of $n$ so that $h_n$ has the same congruence class modulo $D$
where $d_{n-M'}=1$ and $d'_{n-M'}=0$.  (again, if $\beta_{M,i}$ is
negative, we reverse the values of $d$ and $d'$ and use the absolute
value).

By Lemmas \ref{copycolor} and \ref{totaldistancelem} these copies
have the properties that we want.
\end{proof}

\begin{lem}\label{mulreclem}  Given $a,b\in\Gamma$, with some $n$
where $F(n)=a,F(n+1)=b$, and given $k\in \N$ then there exists an
interval $I$ of natural numbers, and functions $f_0,f_1,\ldots,f_k:I
\rightarrow \N_0$ so that:

\begin{enumerate}
\item $0\leq f_i(l) < \gamma_{F(l)}$
\item For every $1\leq i\leq k$, $0\leq x\leq \gamma_a-1$ and $0\leq
y \leq \gamma_b-2$, there exists $n'\in I$ so that
$F(n')=a,F(n'+1)=b$ and the values on $(n',n'+1)$ of $f_0$ and $f_i$
are $(x,y)$ and $(x+1,y)$ (or $(0,y+1)$ if $x=\gamma_a-1$)
respectively.
\item When the $f_i$ are treated as indices for copies of a column of generation
$\min (I)$, these copies are all in the same column and consecutive
copies are separated by the same power of $T$.
\end{enumerate}

\end{lem}

The significance of Lemma \ref{mulreclem} is that it allows us to produce several equally spaced copies of a given interval.  Furthermore the second condition states that the indices of these copies will be rich enough that modifying them will allow us to make crucial adjustments.

\begin{proof}
Find a sequence of consecutive values of $F$ of the form
$e_1,e_2,\ldots,e_w,a,b$ with $\prod_{i=1}^w \gamma_{e_i}>k$. Extend
this to a sequence of the form
$$e_1,e_2,\ldots,e_w,a,b,d_1,\ldots,d_z,e_1,e_2,\ldots,e_w.$$
Find an interval $I$ so that $F$ applied to $I$ yields 
$\left(\prod_{i=1}^n \gamma_{e_i}\right) \gamma_a \gamma_b \left(\prod_{i=1}^l
\gamma_{d_i}\right)-1$ non-intersecting copies of the above sequence.  We
will make our $f_i$ all be 0 off of these subsequences.

We define $f_0$ on these subintervals so that it takes every
possible set of values on the first $w+z+2$ entries (as limited by
property 1) except for all 0's.  On each such block we let $f_0$
take values on the two instances of $e_i$ that add up to
$\gamma_{e_i}-1$.

We think of the values of an $f_i$ on such a block as an appropriate
radix representation (leftmost digit least significant) of a natural
number.  We inductively define $f_{i+1}$ to represent the number one
larger.  In particular if on some such block $f_i$ takes the values
$\gamma_{1}-1,\ldots,\gamma_{s}-1,v_{s+1},\ldots,v_{2w+z+2}$ where
$\gamma_j$ is the appropriate $\gamma$ for the $j^{th}$ term and
$v_{s+1}<\gamma_{s+1}-1$, then on this block $f_{i+1}$ takes values
$0,\ldots,0,v_{s+1}+1,v_{s+2},\ldots,v_{2w+z+2}$.

We note that Property 1 is clearly satisfied.  Property 2 is
satisfied because if we consider the blocks on which $f_0$ has
values $\gamma_{e_1}-1,\ldots,\gamma_{e_w}-1,x,y$, then $f_0$ and
$f_i$ have the appropriate values on the $a,b$ terms.  (Using the
fact that $\prod_{i=1}^n \gamma_{e_i}>i.$) 

We note that by Lemma \ref{consecdistancelem} that the difference in
heights of the consecutive copies indexed by the $f_i$ is a fixed
sum of $h_N$ corresponding to the beginnings of blocks, plus a
correction term based on changes in the number of $N$ for which
$f_i(N)$ and $F(N)$ have particular given values.  Combining this
with Lemma \ref{copycolor} we need only show that the number of such
$N$ remains constant.

For each $f_i$ and each block we associate the three numbers
corresponding to the natural numbers given by the appropriate radix
representations
$f_i(n_0+1),\ldots,f_i(n_0+w)$, and 
$f_i(n_0+w+1),\ldots,f_0(w+z+2)$ and
$f_i(n_0+w+z+3),\ldots,f_0(2w+z+3)$, where $n_0+1$ is the beginning
of the subinterval.  It suffices to show that the multiplicities
with which numbers show up in either of the first and third places
remains constant, and that the multiplicities with which numbers
show up in the second place remains constant.

In $f_i$ the first and second places take all possible values except
for $i,0$.  Now if $M_1,M_2$ are one more than the maximum possible
values in the first and second places, then $f_i$ and $f_0$ agree in
the third place except when the value of $f_0$ in the first two are
$M_1-1-l,M_2-1$ with $0\leq l< i$.  In that case we have carry over
to the third place and there $f_i$ has the value of $l+1$ instead of
$l$.  So in the third place, $f_i$ and an extra $i$ and one fewer
$0$.  This completes our proof.

\end{proof}

We can now prove Theorem \ref{multrecthm}.
\begin{proof}
Let $A$ be a set of positive measure.  Let $k$ be an integer.  From
Lemma \ref{mulreclem} we can see that there is an $\epsilon>0$ and
columns of arbitrarily high generation so that these columns have
$k+1$ copies whose size is more than $\epsilon$ of that of the
original, so that these copies are in the same column and
consecutive copies are separated by the same amount.  Take a level of such a generation that is more than
$\left(1-\frac{\epsilon}{k+1} \right)$-full of $A$.  
Then the copies
of this level in those copies of its column are each more than
$\left(\frac{k}{k+1}\right)$-full of $A$ and have the property that
for some $n$, 
$T^{in}$ of the bottom level is another one of the
levels for $1\leq i\leq k$.  Hence for this $n$, $\mu(A\cap
T^n(A)\cap\cdots\cap T^{kn} (A))>0$.
\end{proof}

We now prove Theorem \ref{pwmthm}.
\begin{proof}
We begin with the levels given to us by Lemma \ref{multdislem}.  We
wish to show that there is an $\epsilon>0$ so that for $m$ pairs of
levels of arbitrarily high generation with the same separation
between corresponding levels as we have between $I_i$ and $J_i$, we
can find copies of these levels of size more than $\epsilon$ times
that of the original so that corresponding copies are in the same
column, and so that the difference in heights between the $i^{th}$
pair of copies is proportional to $k_i$.  This would prove our
Theorem with a simple application of double approximation.

Notice that the intersection of the span of the set in Lemma
\ref{technicallem} with $\Z\times G \times \{ 0\}$ is the span of
$(c_{M,i},0)$, since $t_{N,i}-t_{N,j} = c_{N-1,i}-c_{N-1,j}$.
Therefore, $(D,0)$ is in the span of the $c_{M,i}$.
Hence if $T^{Dd_i}(I_i)=J_i$ then we can write 
$$
(Dd_i,0)=\sum_{M,j} \alpha_{M,j,i} c_{M,j}
$$

Using Lemma \ref{mulreclem} we can find an $\epsilon>0$ so that for
arbitrarily large generations of columns, we can find $m$ pairs of
copies so that each pair of copies is in the same column and
separated by a power of $T$ proportional to $k_i$.  Furthermore, we
can make these copies more than $\epsilon$ the size of the original.
Lastly, we can guarantee that for any $M,j$ with there are at least
$|\alpha_{M,j,i}|$ integers $n$ so that $F(n)=F(M)$, $F(n+1)=F(M+1)$
and so that in the indexing of the first and  second copy in the
$i^{th}$ pair, the index of the copy at the digits corresponding to
$n$ and $n+1$ are either $0$,$j$ and $1$,$j$ or $\gamma_n -1$,$j$
and $0$,$j+1$ respectively.  If we change $|\alpha_{M,j,i}|$ of one
of these types to the other, then we keep these pairs of copies in
the same column, but alter their relative height difference by
$Dd_i$.  This provides what we need for the double approximation.
\end{proof}

\begin{remark}
Note that the proof of Theorem \ref{pwmthm} implies that $T$ is
power weakly mixing if and only if it is totally ergodic.  Note also
that to check condition 2, it is sufficient to first check to see if
a $D$ from Lemma \ref{technicallem} exists, and if one does, to check
condition 2 modulo $D$ for a particular $N$ so that infinitely often
$F(M)=F(N)$ and $h_N \equiv h_M \pmod{D}$.  This reduces checking
condition 2 to a finite computation.

Notice also that if $\emph{\textrm{Im} (F)}=\{
(n,s_0,\ldots,s_{n-1},0,g_1,\ldots,g_{n-1})\}$, then condition 2 can
be written in the simple form that $(1,0)$ is in the span of
$\{(s_i+h_N,g_{i+1}-g_{i}),(s_{n-1},-g_{n-1})\}$ for all $N$.
\end{remark}

\section{Examples}\label{S:examples}
These first few will be where $\Gamma$ has a single element as in the
last remark.

Consider for example the  Chac\'on-$m$ transformation for $m\geq 2$.
Following N. Friedman, a Chac\'on-$m$ transformation is a rank-one
transformation where column $C_{n+1}$ is obtained from column
$C_n$ but cutting each level of $C_n$ into $m$ sub-levels, stacking
from left to right and placing a spacer on top of the last level;  see
Section~\ref{S:pci} for more details on rank-one constructions.
We can
 define a Chac\'on-$m$ transformation by letting $G=\{0\}$, $n=m$ and $s_i=0$ for $0\leq i \leq
n-2$, $s_{n-1}=1$, and $g_i=0$ for $0\leq i \leq n-1$.  Clearly, $\{
g_i\}$ generates $G$, so condition 1 is satisfied.  Since $(1,0) =
(s_{n-1},-g_{n-1})$ it is in the span of
$\{(s_i+h_N,g_{i+1}-g_i),(s_{n-1},-g_{n-1}) \}$.  Therefore, it is
power weakly mixing.

Consider the transformation defined by $G=\{0\}$, $n=3$ and the
sequence $(1,1,0,0,0,0)$.  Condition 2 states that $(1,0)$ is in the
span of $\{(1+h_N,0),(1+h_N,0),(0,0) \}= (1+h_N)\Z \times \{ 0\}$.
Which does not hold for any, $h_N$.  Therefore this transformation
is not power weakly mixing (in fact it is not $T^2$ ergodic).

Consider the transformation, $T$ defined by the group $G=\Z$, $n=5$
and the set $(0,0,0,1,0,0,1,0,0,0)$.  It satisfies condition 1 since
$g_1=1$ generates $G$.  Condition 2 states that $(1,0)\in
\textrm{span}\{ (h_N,1),(h_N,-1),(h_N,0),(1+h_N,0),(0,0)\}$.  This
clearly holds since $(1,0)=(1+h_N,0)-(h_N,0)$.  Therefore, $T$ is a
power weakly mixing, infinite-measure preserving transformation.

Lastly consider the group $G$ to be any countably generated abelian
group with generators $e_i$ for $i\in \No$.  If $n\in\N$, let
$e(n)=k$ where $k$ is the largest power of 2 such that $2^k$ divides
$n$.  Let $F(n)=(4,0,0,0,1,0,0,e_{e(n+1)},0)$.  $F$ clearly
satisfies the necessary condition.  Notice that $T$ is a
$G$-extension of Chac\'on-$4$.  Condition 1 is clearly satisfied.
Condition 2 is satisfied since
\begin{align*}
(1,0) = (& (s(1,0)+s(0,3)-s(0,0)), (g(1,1)-g(1,0)-g(0,3)-g(0,1))).
\end{align*}
Therefore, $T$ is power weakly mixing.


\section{Non-totally ergodic $2$-point extension}


As an example of the above we analyze what happens in the particular case where $G=\Z/2\Z$ and $F(n)=(2,0,1,0,1)$.

Consider the two-point extension $T$ of the Chacon-$2$
transformation formed as follows.  Begin with two intervals of equal
size -- call them columns $C_{0,0}$ and $C_{0,1}$.  These will be
known as the generation zero columns.  To define the generation
$n+1$ columns, cut each of the generation $n$ columns in half,
stacking the right half of $C_{n,1}$ over the left half of $C_{n,0}$
and vice versa. Then add a spacer to the top of the two columns thus
formed to yield the generation $n+1$ columns.  The transformation $T$ is defined to map each point to the
point directly above it. To see that $T$ is indeed a two-point
extension of Chacon-$2$, associate each point in $C_{n,0}$ with the
corresponding point in $C_{n,1}$ for all $n$. The resulting space
and transformation are exactly Chacon-$2$. It is well-known that
Chacon-$2$ is a weakly mixing transformation, therefore totally
ergodic, but we will show that the two-point extension $T$ is not
even $T^2$ ergodic.

Let $I$ and $J$ be the top and middle levels, respectively, in
$C_{1,0}$.  Suppose for some $m$ we have $\mu(T^m(I)\cap J)>0$. Then there must
be some generation in which there is a copy $J'$ of $J$ above a copy
$I'$ of $I$ by a distance $m$ levels.  That is, $T^m(I') = J'$, and
we will write $d(I',J')=m$. We will show that this cannot be the case if $m$ is even.

We prove by induction on $n$ that for all $n$, any two copies of $I$ in one of the generation $n$
columns must be an even distance apart.  This is vacuously true for
$n = 1,2$.  Assume it is true for generation $n$.  Label the left
halves of the top copies of $I$ in $C_{n,0}$ and $C_{n,1}$ as $I_0$
and $I_1$, respectively.  Label the right halves of the bottom
copies of $I$ as $I_2$ and $I_3$, respectively.  To prove the claim
we must show that the distances from $I_0$ to $I_3$ and from $I_1$
to $I_2$ are both even.  For convenience label the right halves of
the bottom levels of $C_{n,0}$ and $C_{n,1}$ by $K_0$ and $K_1$
respectively.  Then $d(K_0,I_2)=2$ and $d(K_1,I_3)=5$
since this is true for $n=2$ and since the bottoms of columns are preserved through later generations. The distance
from $I_0$ to $K_1$ is given by:
\[
  d(I_0,K_1)= \begin{cases} n + 3\ \ \ & n \textrm{ even}\\
                             n\ \ \ & n \textrm{ odd}
                \end{cases}
\]
Similarly, we have
\[
 d(I_1,K_0)= \begin{cases} n\ \ \ & n \textrm{ even}\\
                            n + 3\ \ \ & n \textrm{ odd}
                \end{cases}
\]
These are true because they hold for $n=2$ and by induction on $n$.  The fundamental idea is that the top of $C_{n,1}$ looks like the top of $C_{n-1,0}$ with an extra spacer added on top.  From this the above statements are easily shown by induction.

Both $d(I_0,I_3)=d(I_0,K_1) + d(K_1,I_3)$ and
$d(I_1,I_2)=d(I_1,K_0)+ d(K_0,I_2)$ must then be even, independent
of $n$. By induction, $\mu(T^m(I)\cap I)>0$ implies $m$ is
even.  Since each copy of $J$ lies directly below a copy of $I$,
this means that $\mu(T^m(I)\cap J)>0$ implies $m$ is odd.
Therefore $\mu(T^{2m}(I)\cap J) = 0$ for all $m$, so $T^2$ is not
ergodic.



\section{Conservativity and Recurrence on a Sufficient Class}

In this section all transformations are assumed to be infinite measure preserving, and not necessarily invertible.  If for any measurable set $A$ we have
$\mu(A \setminus \bigcup_{i=1}^{\infty}T^{-i}A) = 0$ then $T$ is
said to be \textbf{recurrent}.  For sets of finite measure this condition is equivalent to $\mu(A \cap \bigcup_{i=1}^{\infty}T^{-i}A) = \mu(A)$.
A class $\mathcal{C}$ of subsets of
$X$ is called a \textbf{sufficient class} if it satisfies the
following approximation property for all measurable $A \subset X$:
\[
\mu(A)=\inf\left\{ \sum_{j=1}^{\infty}\mu(I_j):\{I_j\}\textrm{ cover
A and }I_j\in\mathcal{C}\right\}
\]
A transformation is said to be \textbf{conservative on
$\mathcal{C}$} or \textbf{recurrent on $\mathcal{C}$} if the
condition for conservativity or recurrence holds for all $I\in\mathcal{C}$ of positive measure, but not necessarily for all measurable sets.
While conservativity and recurrence are known to be equivalent, we show in this section that
conservativity and recurrence on a sufficient class $\mathcal{C}$
are not equivalent.  In particular recurrence on a sufficient class implies recurrence, but the same is not true for conservativity.

Consider the following infinite measure preserving transformation
$T: \Rm\rightarrow\Rm$ which is conservative on the sufficient class
$\mathcal{C}=\{I: I\textrm{ is a finite open interval}\}$.  It is
well known that there exist sets $K\subset [0,1)$ and
$K^c=[0,1)\setminus K$ of positive measure such that $\mu(I\cap
K)>0$ if and only if $\mu(I\cap K^c)>0$ for all $I\in \mathcal{C}$.
Define $T$ by:
\[
 T(x) = \begin{cases} x\ \ \ & x\mod 1\in K\\
                        x+1\ \ \ & x\mod 1 \in K^c
        \end{cases}
\]
Then $T$ is conservative on $\mathcal{C}$, but $\mu(T^{-n}(K^c)\cap
K^c)=0$ for all nonzero $n$ so $T$ is not conservative.  Note
however that $T$ is not recurrent on $\mathcal{C}$.

\begin{prop}
Let $(X,\mathcal{B},\mu)$ be a measure space such that any set of infinite measure has a subset of finite but positive measure, e.g. a $\sigma$-finite measure space, and let $T$ be an infinite measure-preserving transformation.  If $T$ is recurrent on
a sufficient class, then $T$ is recurrent.
\end {prop}

\begin{proof}
Let $\mathcal{C}$ be a sufficient class, and suppose that $T$ is not
recurrent.  Then $T$ is not conservative, so there exists a set $A$
of positive measure such that $\mu(A\cap T^{-n}(A)) = 0$ for all $n
> 0$.  Perhaps taking a subset, we may assume 
that $A$ has finite measure.  We can then find a set $I \in \mathcal{C}$ such that $\mu(A\cap I)
> \frac{1}{2}\mu(I)$.  Note that if a subset of $I
\setminus A$ with positive measure is mapped into $A$ by $T^{-i}$ for some $i$ it will never be mapped into $A$ for any $j>i$ (by hypothesis, any subset of $A$ of positive measure is never mapped into $A$ under iteration by $T^{-1}$), so we have
\begin{align*}
\mu\left(I \cap \bigcup_{i=1}^{\infty}T^{-i}(I)\right) & =
\mu\left((I\setminus A) \cap \bigcup_{i=1}^{\infty}T^{-i}(I)\right)
+ \mu\left((I\cap A) \cap \bigcup_{i=1}^{\infty}T^{-i}(I\setminus A)\right) \\
& \leq \mu\left(I \setminus A\right) +  \mu\left(I \setminus A\right)\\ & <
\frac{1}{2}\mu(I) + \frac{1}{2}\mu(I) = \mu(I)
\end{align*}
Thus $T$ is not recurrent on $\mathcal{C}$.
\end{proof}

To see that a regularity condition on the space $(X,\mathcal{B},\mu)$ is necessary, let $X=\Rm$ and $\mathcal{B} = 2^\Rm$.  Define $\mu(A)$ to be zero if $A$ is finite or countable and infinity otherwise.  The collection $\mathcal{C}$ of all singletons in $\Rm$ along with the set $\Rm$ itself is a sufficient class for this space.  Any bijective map on $X$ is measure-preserving and recurrent on $\mathcal{C}$, but in general it need not be recurrent.


\section{Power Conservative Index}  \label{S:pci}
	In this section we  obtain a condition for rank-one transformations that
implies power conservative index.  A notion that has been used to study conservativity of products is that of positive type.  A transformation $T$ is of {\bf positive type} if
 $\limsup_{n\to\infty} \mu(T^n(A)\cap A)>0$   for all sets $A$ of
positive measure.  Clearly, if
$T$ is of positive type, then it is conservative.  It was  shown in \cite{AN00} that if $T$ is of positive
type,  then for each positive integer $d$, the cartesian product of $d$ copies of  $T$ is of positive type, so positive type implies infinite conservative  index.   But it  is easily verified that the transformation $T_1$ of \cite{AFS01}   is  of positive type
 but as already mentioned $T_1 \times T_1^2$ is not conservative, so $T_1$ is  not power conservative,  showing that positive type does not imply power conservative index. Our condition can be  used to show show that some infinite measure-preserving staircases have power conservative index.

	First, we introduce some  notation for  constructing
	 measure-preserving  rank-one transformations.
We start with a column $C_0$, which is a unit interval. Let ${r_n}$ be a
sequence of integers with $r_n\geq 2$. At each stage $n$, we have a column
$C_n$ that consists of $h_n$ intervals ($h_n$ denotes the height of
column $C_n$, which is definite to be the number of intervals in the column). We denote
the intervals in $C_n$ by $I_{n,0}, I_{n,1},..., I_{n,h_n-1}$, where $I_{n,0}$ is the
interval at the lowest level. A column determines a map on all
of its levels but the top, where
each interval is mapped  to the interval above it by the canonical translation.  Column $C_{n+1}$
is obtained from column $C_n$
by cutting and stacking according to the following procedure. Cut
all intervals of column $C_n$ into $r_n$ subintervals of the same measure
 to form subcolumns, $C_n^{[0]},
C_n^{[1]},..., C_n^{[r_{n}-1]}$.
Then we may put spacers on top of the subcolumns. Let
$\{s_{n,i}\}_{i=0}^{r_n-1}$ be a doubly indexed sequence of non-negative
integers. The sequence $s_{n,0}, s_{n,1},..., s_{n,r_n-1}$
specifies the number of spacers on each respective subcolumn. Then  stack
subcolumns on top of one another, with each subcolumn going underneath
its adjacent subcolumn to the right, so that the rightmost subcolumn goes
on the very top. This cutting and stacking procedure
obtains  column $C_{n+1}$. This defines a sequence of columns ${C_n}$
and as the width of the column approaches $0$, it defines a measure-preserving
transformation on a finite or an infinite interval.

We prove the following Theorem:
\begin{thm}
\label{T: 2.1}
Let $T$ be a rank 1 transformation with cut sizes $\{r_n\}$.  If for all $d>0$,
$$
\liminf_{n\rightarrow\infty} \frac{h_n^{d-1}}{\prod_{i=0}^{n-1}r_i^d} = 0
$$
then $T$ is power conservative index.
\end{thm}

We do this by proving the stronger theorem
\begin{thm}
\label{stronger conservativity theorem}
Let $T$ be a rank 1 transformation with cut sizes $\{r_n\}$.  Let $\{k_i\}_{i=1}^d$ be positive integers.  If
$$
\liminf_{n\rightarrow\infty} \frac{h_n^{d-1}}{\prod_{i=0}^{n-1}r_i^d} = 0,
$$
then $T^{k_1}\times T^{k_2} \times \ldots \times T^{k_d}$ is conservative.
\end{thm}

Fix the $k_i$.  Let $S=T^{k_1}\times T^{k_2} \times \ldots \times T^{k_d}$.  For a column $C_n$ define an equivalence relation $\sim_n$ on $d$-fold products of levels so that ${I} \sim_n {J}$ if and only if for some integer $N$, $S^N({I}) = {J}$.

\begin{lem}
\label{A size lem}
Suppose that $A\subset X^d$ is a set where $\mu(S^n(A)\cap A) = 0$ for all $n\neq 0$.  Let $L$ be the $\sim_n$ equivalence class of products of levels equivalent to $I$.  Then
$$
\mu\left(A\cap\left(\bigcup_{J\in L} J\right) \right) \leq \mu(I).
$$
\end{lem}
\begin{proof}
First modify $A$ by a set of measure 0 so that $S^n(A)\cap A = \emptyset$ for all $n\neq 0$.
Note that $\bigcup_{J\in L} J$ is a subset of $\bigcup_{N\in\Z} S^N(I)$.  In fact for some subset $P\subset \Z$, we can write 
$$
\bigcup_{J\in L} J = \coprod_{N\in P} S^N(I).
$$
We therefore, think of this set as $P\times I$.  Let $\chi_A$ be the characteristic function of $A$ on $P\times I$.  We note that since $S^n(A)\cap A = \emptyset$ for $n\neq 0$, that $\chi_A(a,x)$ and $\chi_A(b,x)$ cannot both be 1 for $a\neq b$.  Hence
$$
\int_P \chi_A(n,x) dn = \sum_{n\in P} \chi_A(n,x) \leq 1.
$$
So by changing the order of integration, we get that
\begin{align*}
\mu\left(A\cap\left(\bigcup_{J\in L} J\right) \right) & = \int_{P\times I} \chi_A(n,x) d\mu\\
& = \int_I \int_P \chi_A(n,x) dn dx\\
& \leq \int_I dx \\
& = \mu(I).
\end{align*}
\end{proof}

\begin{lem}
\label{number of equiv classes}
The number of equivalence classes of $\sim_n$ is at most 
$$\left(\sum_{i=1}^d k_i \right)h_n^{d-1}.$$
\end{lem}
\begin{proof}
To each product of levels, we can associate a $d$-tuple of integers in the range $[1,h_n]$ representing the heights of the levels.  It is clear that the product associated with $\{a_i\}$ and the product associated with $\{a_i+k_i\}$ are equivalent.  Each equivalence class has at least one element with $\sum_{i=1}^d a_i$ minimal.  We will bound the number of such sequences.

Clearly, for $\sum_{i=1}^d a_i$ to be a minimal representative of an equivalence class, $\{a_i-k_i\}$ cannot be a valid sequence.  Therefore $a_i\leq k_i$ for some $i$.  The number of such sequences with $a_i\leq k_i$ for a given $i$ is at most $k_i h_n^{d-1}$.  Summing over $i$ gives our result.
\end{proof}

We are now prepared to prove Theorem \ref{stronger conservativity theorem}.

\begin{proof}
Suppose that 
$$
\liminf_{n\rightarrow\infty} \frac{h_n^{d-1}}{\prod_{i=0}^{n-1}r_i^d} = 0,
$$
and that $A\subset X^d$ satisfies $\mu\left(S^n(A)\cap A\right) = 0$ for all $n\neq 0$.  We wish to bound $\mu\left(A\cap C_n^d\right)$.  Lemma \ref{A size lem} says that the intersection of $A$ with any equivalence class of levels is at most the size of a level, or $$\left(\prod_{i=0}^{n-1}r_i^d\right)^{-1}$$.  Lemma \ref{number of equiv classes} says that there are at most $K h_n^{d-1}$ equivalence classes where $K=\sum_{i=1}^d k_i$.  Therefore,
$$
\mu\left(A\cap C_n^d\right) \leq \frac{K h_n^{d-1}}{\prod_{i=0}^{n-1}r_i^d}.
$$
Hence
$$
\liminf_{n\rightarrow\infty} \mu\left(A\cap C_n^d\right) = 0.
$$
Since $C_n^d$ exhausts $X^d$, this implies that $\mu(A)=0.$
\end{proof}

\begin{cor}
\label{C: 9.9}
There exists an infinite measure-preserving pure staircase
transformation $T$ such
that $T$ has power conservative index.
\end{cor}
\begin{proof}
Let $T$ be the classical staircase transformation with $r_m =
2^{2^m}$.  Direct computation shows that $h_m\leq \frac{m+1}{2} r_m$ holds for $m = 0$, and we will prove that it holds for all $m$ by induction.  If this holds for $m$ we can apply the definition of $h_{m+1}$ for a staircase transformation to yield
\begin{equation*}
	h_{m+1} = h_m r_m + \frac{1}{2}r_m(r_m-1) \leq \frac{m+1}{2} r_m^2 + \frac{1}{2}r_m^2 = \frac{m+2}{2}r_m^2 = \frac{m+2}{2} r_{m+1}.
\end{equation*}
Then we can bound
\begin{equation*}
	\frac{h_m^{d-1}}{(\prod_{i=0}^{m-1}r_i)^d} = \frac{h_m^{d-1}}{\left(\frac{r_m}{2}\right)^d} = \frac{1}{h_m}\left(\frac{2h_m}{r_m}\right)^d \leq \frac {(m+1)^d}{h_m} \leq \frac{(m+1)^d}{r_m}
\end{equation*}
which approaches $0$ as $m\rightarrow\infty$ for all $d\geq 1$.  Hence, this transformation $T$ satisfies the
condition of Theorem \ref{T: 2.1}.
\end{proof}

\begin{cor}
Let $T$ be a rank-one  transformation with $r_n = 2^{2^n}$ and such that
\[
s_{n,i} = \begin{cases} 2i \ \  & \text{when}\ 2 | i \\
                        0 \ \  & \text{otherwise.}
          \end{cases}
\]
Then $T$ is an infinite measure-preserving
transformation that is power conservative but such that $T^2$ is not ergodic.
\end{cor}

\begin{proof}
 An argument similar to that in Corollary \ref{C: 9.9} shows that $T$ has power conservative index.  Let $I_1$ and $I_2$ be two levels in some column $C_m$ such that
 $\alpha (I_1) -\alpha (I_2)$ is a positive odd number.  Now consider sublevels of $I_1$ and $I_2$, denoted by $I_1^\prime$ and $I_2^\prime$ respectively,
 in some column $C_n$ for $n>m$. From the construction, $\alpha (I^\prime_1) -\alpha (I^\prime_2)$ will
 always be an odd number. It follows that there does not exist an integer $t$ such that $T^{2t}(I_1^\prime) = I_2^\prime$.  Hence, $T^2$ is not ergodic.
 \end{proof}


\begin{thebibliography}{99}

\bibitem{A98}  T. Adams.
{Smorodinsky's conjecture on rank-one mixing},    {\it  Proc. 
 Amer. Math. Soc. }  {\bf 126(3)}  (1998), 739-744.

\bibitem{AN00}
        J. Aaronson and H. Nakada. 
        { Multiple recurrence of Markov
     shifts and other infinite measure preserving transformations. {\it 
      Isr.  J.  Math.}  117, (2000), 285-310.

    \bibitem{AFS97}
    T. Adams, N. Friedman, and C.E. Silva. {Rank-one weak mixing
    for nonsingular transformations}. {\it  Isr. J.  of Math.}  102
(1997),
    269-281.

    \bibitem{AFS01}
    T. Adams, N. Friedman, and C.E. Silva.  {Rank one power
weakly
    mixing nonsingular transformations},{\it  Ergodic Theory \& Dynam.
    Sys.}  21 (2001),  1321-1332.

    \bibitem{BFMS}
    A. Bowles, L. Fidkowski, A. Marinello, and C.E. Silva. 
    {Double ergodicity of nonsingular transformations and
    infinite measure-preserving staircase transformations}. {\it Illinois
    J.  of Math.}  45, ( 2001),  999--1019.

    \bibitem{D}
A. Danilenko.
 { Funny rank-one weak mixing for nonsingular abelian actions. }
{\it Israel J. Math.} 121 (2001), 29--54.

    \bibitem{DGMS}
    S. Day, B. Grivna, E. McCartney, and C.E. Silva. {Power
    Weakly Mixing Infinite Transformations}. {\it  New York Journal of
    Math.}  { 5} (1999), 17-24.


    \bibitem{EHH98} S. Eigen, A. Hajian, K. Halverson, {
    Multiple recurrence and Infinite Measure Preserving Odometers},
    {\it Israel J. Math.} 108 (1998), 37-44.


    \bibitem{hF81}
    H. Furstenberg. {\it Recurrence in ergodic theory and
combinatorial number
    theory}. 
    Princeton Univ. Press, Princeton, N.J., 1981.

 \bibitem{G03}
  K. Gruher, F. Hines, D. Patel, C. E. Silva and R. Waelder.
{Power weak mixing does not imply multiple recurrence in infinite
 measure and other counterexamples}. {\it  New York J. Math.} 
     {9},
       {2003}, 1--22.

\bibitem{inoue}    K. Inoue.  {Isometric extensions and multiple recurrence
of infinite measure preserving systems} }{\it  Israel J.
Math.} 140 (2004) 245--252.

\bibitem{Mey} T. Meyerovitch.  Extensions and multiple recurrence of infinite measure  
preserving systems, preprint. ArXiv: 
http://arxiv.org/abs/math/0703914.

\bibitem {Ro} E. A. Robinson. {{Ergodic properties that lift to compact group extensions}},
  {\it  Proc. Amer. Math. Soc.},
    {102},
      {1988},
61-- 67.

\bibitem{Ru}D. J. Rudolph. { {$k$-fold mixing lifts to weakly mixing isometric extensions}}.
{\it Ergodic Theory \& Dynam. Systems}  {\bf 5} (1985), no.~3, 445--447

\bibitem{Si08} C. E. Silva. {\it Invitation to Ergodic Theory}. Student Math. Library, Vol. 42, Amer. 
Math. Soc., 2008.

\end{thebibliography}
\end{document}